\documentclass[preprint,12pt]{elsarticle}

\usepackage{graphicx,color}

\usepackage{amssymb}
\usepackage{amsthm}

\usepackage{amsmath,mathrsfs,bm,subfigure,verbatim}

\journal{Journal of Computational and Applied Mathematics}

\newtheorem{thm}{Theorem}
\newtheorem{lem}[thm]{Lemma}
\newtheorem{prop}[thm]{Proposition}
\newtheorem{cor}[thm]{Corollary}
\newdefinition{defn}{Definition}
\newdefinition{rmk}{Remark}
\newdefinition{alg}{Algorithm}
\newdefinition{exmp}{Example}
\newproof{pf}{Proof}

\begin{document}

\begin{frontmatter}




\title{A Bivariate Preprocessing Paradigm for Buchberger-M\"{o}ller Algorithm}
\tnotetext[fund]{This work was supported in part by the National
Grand Fundamental Research 973 Program of China (No. 2004CB318000).}


\author{Xiaoying Wang, Shugong Zhang}
\author{Tian Dong\corref{cor1}}
\cortext[cor1]{Corresponding author} \ead{dongtian@jlu.edu.cn}
\address{School of Mathemathics,
Key Lab. of Symbolic Computation and Knowledge Engineering
\textup{(}Ministry of Education\textup{)}, Jilin University,
Changchun 130012, PR China}

\begin{abstract}
For the last almost three decades, since the famous Buchberger-M\"{o}ller(BM) algorithm
emerged, there has been wide interest in vanishing ideals of points and associated
interpolation polynomials. Our paradigm is based on the theory of bivariate polynomial
interpolation on cartesian point sets that gives us related degree reducing interpolation
monomial and Newton bases directly. Since the bases are involved in the computation
process as well as contained in the final output of BM algorithm, our paradigm obviously
simplifies the computation and accelerates the BM process. The experiments show that the
paradigm is best suited for the computation over finite prime fields that have many
applications.

\end{abstract}


\end{frontmatter}


\section{Introduction}
\label{sec:int}

For an arbitrary field $\mathbb{F}$, we let  $\mathbb{F}_q$ a finite prime field of size
$q$ and $\Pi^d:=\mathbb{F}[x_1, \ldots, x_d]$ the $d$-variate polynomial ring over
$\mathbb{F}$. Given a preassigned set of distinct affine points $\Xi\subset
\mathbb{F}^d$, it is well-known that the set of all polynomials in $\Pi^d$ vanishing at
$\Xi$ constitutes a radical zero-dimensional ideal, denoted by $\mathcal{I}(\Xi)$, which
is called the \emph{vanishing ideal} of $\Xi$.

Recent years, there has been considerable interest in vanishing ideals of points in many
branches of mathematics such as algebraic geometry\cite{CLO2005}, multivariate
interpolation\cite{Sau2006, dBo2005}, coding theory\cite{Sak1998, Sal2009},
statistics\cite{Rob1998}, and even computational molecular biology\cite{LS2004, JS2006}.
As is well known, the most significant milestone of the computation of vanishing ideals
is the algorithm presented in \cite{MB1982} by Hans Michael M{\"{o}}ller and Bruno
Buchberger known as Buchberger-M\"{o}ller algorithm(BM algorithm for short). For any
point set $\Xi\subset \mathbb{F}^d$ and fixed term order $\prec$, BM algorithm yields the
reduced Gr\"{o}bner basis for $\mathcal{I}(\Xi)$ w.r.t. $\prec$ and a $\prec$-degree
reducing interpolation Newton basis for $d$-variate Lagrange interpolation on $\Xi$. The
algorithm also produces the Gr\"{o}bner \'{e}scalier of $\mathcal{I}(\Xi)$ w.r.t. $\prec$
as a byproduct. Afterwards, in 1993, BM algorithm was applied in \cite{FGLM1993} in order
to solve the renowned FGLM-problem. In the same year, \cite{MMM1993} merged BM and FGLM
algorithms into four variations that can solve more general zero-dimensional ideals
therefore related ideal interpolation problems \cite{dBo2005}. The algorithms are
referred as MMM algorithms.

Although very important, BM algorithm (and MMM algorithms) has a very poor complexity
that limits its applications. In this decade, many authors proposed new algorithms that
can reduce the complexity but mostly suitable for special cases. \cite{ABKR2000}
presented a modular version of BM algorithm that is best suited to the computation over
$\mathbb{Q}$. \cite{CM1995, GRS2003, FRR2006} presented algorithms for obtaining, with
relatively little effort, the Gr\"{o}bner \'{e}scalier of a vanishing ideal w.r.t. the
(inverse) lexicographic order that can lead to an interpolation Newton basis or the
reduced Gr\"{o}bner basis for the vanishing ideal after solving a linear system.

For a fixed point set $\Xi$ in $\mathbb{F}^d$ and a term order $\prec$, it is well known
that there are \emph{two} factors that determine the Gr\"{o}bner \'{e}scalier of
$\mathcal{I}(\Xi)$ w.r.t. $\prec$ thereby the reduced Gr\"{o}bner basis for
$\mathcal{I}(\Xi)$ and related  degree reducing interpolation Newton bases (up to
coefficients). One is apparently the cardinal of $\Xi$. It is the unique determinate
factor in univariate cases. Another one is the geometry (the distribution of the points)
of $\Xi$ that is dominating in multivariate cases but not taken into consideration by BM
and MMM algorithms. Recent years, \cite{Sau2004, Cra2004, CDZ2006} studied multivariate
Lagrange interpolation on a special kind of point sets, cartesian point sets (aka lower
point sets), and constructed the associated Gr\"{o}bner \'{e}scalier and degree reducing
interpolation Newton bases theoretically. We know from \cite{MB1982, MMM1993} that, for a
cartesian subset of $\Xi$ (it always exists!), certain associated degree reducing
interpolation Newton basis forms part of the output of BM algorithm w.r.t. some
reordering of $\Xi$. Therefore, finding a large enough cartesian subset of $\Xi$ with
little enough effort will reduce the complexity of BM algorithm.

Following this idea, the paper proposes a preprocessing paradigm for BM algorithm with
the organization as follows. The next section is devoted as a preparation for the paper.
And then, the main results of us are presented in two sections. Section \ref{special}
will pursue the paradigm for two special term orders while Section \ref{general} will set
forth our solution for other more general cases. In the last section, Section
\ref{Application}, some implementation issues and experimental results will be
illustrated.

\section{Preliminary}\label{pre}

In this section, we will introduce some notation and recall some basic facts for the
reader's convenience. For more details, we refer the reader to \cite{BW1993, Lor1992}.

We let $\mathbb{N}_0$ denote the monoid of nonnegative integers. A \emph{polynomial}
$f\in \Pi^2$ is of the form
$$
f=\sum_{\bm{\alpha}\in \mathbb{N}_0^2} f_{\bm{\alpha}} X^{\bm{\alpha}},
\hspace{0.8cm}\#\{\bm{\alpha}\in \mathbb{N}_0^2: 0\neq f_{\bm{\alpha}}\in \mathbb{F} \} <
\infty,
$$
where \emph{monomial} $X^{\bm{\alpha}}=x^{\alpha_1}y^{\alpha_2}$ with
$\bm{\alpha}=(\alpha_1,\alpha_2)$. The set of bivariate monomials in $\Pi^2$ is denoted
by $\mathbb{T}^2$.

Fix a term order $\prec$ on $\Pi^2$ that may be lexicographical order $\prec_\text{lex}$,
inverse lexicographical order $\prec_\text{inlex}$, or total degree inverse
lexicographical order $\prec_\text{tdinlex}$ etc. For all $f \in \Pi^2$, with $f \neq 0$,
we may write
$$
f = f_{\bm{\gamma_1}} X^{\bm{\gamma_1}} + f_{\bm{\gamma_2}} X^{\bm{\gamma_2}} + \cdots +
f_{\bm{\gamma_r}} X^{\bm{\gamma_r}},
$$
where $0 \neq f_{\bm{\gamma_i}} \in \mathbb{F}, \bm{\gamma_i}\in \mathbb{N}_0^2,
i=1,\ldots,r$, and $X^{\bm{\gamma_1}} \succ X^{\bm{\gamma_2}} \succ \cdots \succ
X^{\bm{\gamma_r}}$. We shall call $\mathrm{LT}(f):=f_{\bm{\gamma_1}} X^{\bm{\gamma_1}}$
the \emph{leading term} and $\mathrm{LM}(f):=X^{\bm{\gamma_1}}$ the \emph{leading
monomial} of $f$. Furthermore, for a non-empty subset $F \subset \Pi^2$, put
$$
\mathrm{LT}(F) := \{ \mathrm{LT}(f) : f \in F\}.
$$
As in \cite{dBo2007}, we define the $\prec\mspace{-8mu}-$\emph{degree} of a polynomial
$f\in \Pi^2$ to be the leading bidegree w.r.t. $\prec$
$$
\delta(f):=\bm{\gamma},\quad X^{\bm{\gamma}}=\mathrm{LM}(f),
$$
with $\delta(0)$ undefined. Further, for any finite dimensional subset $F \subset \Pi^2$,
define
$$
\delta(F):=\max_{f\in F}\delta(f).
$$
Finally, for any $f, g\in \Pi^2$, if $\delta(f) \prec \delta(g)$ then we say that $f$ is
of \emph{lower degree} than $g$ and use the abbreviation
$$
f\prec g:=\delta(f) \prec \delta(g).
$$
In addition, $f\preceq g$ is interpreted as the degree of $f$ is lower than or equal to
that of $g$.

Let $\mathcal{A}$ be a finite subset of $\mathbb{N}_0^2$. $\mathcal{A}$ is called a
\emph{lower} set if, for any $\bm{\alpha}=(\alpha_1,\alpha_2) \in \mathcal{A}$, we always
have
$$
\mathrm{R}(\bm{\alpha}):= \{(\alpha'_1,\alpha'_2)\in \mathbb{N}_0^2: 0\leq \alpha'_i \leq
\alpha_i, i=1,2\} \subset \mathcal{A}.
$$
Especially, $\bm{0} \in \mathcal{A}$. Moreover, we set $m_j=\max_{(h,j)\in \mathcal{A}}h,
0\leq j\leq \nu$, with $\nu=\max_{(0,k)\in \mathcal{A}} k$. Clearly, $\mathcal{A}$ can be
determined uniquely by the ordered $(\nu+1)$-tuple $(m_0,m_1, \ldots,m_{\nu})$ hence
represented as $\mathrm{L}_x(m_0,m_1,\ldots,m_{\nu})$. Swapping the roles of $x$ and $y$,
we can also represent $\mathcal{A}$ as $\mathrm{L}_y(n_0, n_1, \ldots, n_{m_0})$ with
$n_i=\max_{(i,k)\in \mathcal{A}}k, 0\leq i \leq m_0$. It should be noticed that
$\nu=n_0$.

Given a set $\Xi=\{\xi^{(1)}, \ldots, \xi^{(\mu)}\} \subset \mathbb{F}^2$ of $\mu$
distinct points. For prescribed values $f_i \in \mathbb{F}, i=1, \ldots, \mu$, find all
polynomials $p \in \Pi^2$ satisfying
\begin{equation}\label{inter}
p(\xi^{(i)}) = f_i,\quad i=1,\ldots,\mu.
\end{equation}
We call it the problem of  \emph{bivariate Lagrange interpolation}.
Note that in most cases, especially from a numerical point of view,
we are not interested in all such $p$'s but a ``degree reducing"
one, as in the univariate cases.

\begin{defn}\cite{Sau2006}\label{MDIS}
Fix term order $\prec$. We call a subspace $\mathcal{P} \subset \Pi^2$ a \emph{degree
reducing interpolation space} w.r.t. $\prec$ for the bivariate Lagrange interpolation
(\ref{inter}) if

\begin{description}
    \item[\textbf{DR1}.] $\mathcal{P}$ is an \emph{interpolation space}, i.e., for any $f_i \in \mathbb{F},
          i=1, \ldots, \mu$, there is a unique $p \in \mathcal{P}$ such that $p$ satisfies (\ref{inter}). In
           other words, the interpolation problem is \emph{regular} w.r.t. $\mathcal{P}$.

    \item[\textbf{DR2}.] $\mathcal{P}$ is $\prec\mspace{-8mu}-$\emph{reducing}, i.e., when
          $L_{\mathcal{P}}$ denotes the Lagrange projector
          with range $\mathcal{P}$, then the interpolation polynomial
          $$
          L_{\mathcal{P}} q \preceq q, \hspace{5mm}\forall q \in \Pi^2.
          $$
\end{description}

\end{defn}

For interpolation problem (\ref{inter}), a given interpolation space $\mathcal{P}\subset
\Pi^2$ will give rise to an \emph{interpolation scheme} that is referred as $(\Xi,
\mathcal{P})$, cf. \cite{Lor1992}. Since (\ref{inter}) is regular w.r.t. $\mathcal{P}$,
we can also say that $(\Xi, \mathcal{P})$ is regular. Moreover, if $\mathcal{P}$ is
degree reducing w.r.t. $\prec$, a basis $\{p_1, \ldots, p_\mu\}$ for $\mathcal{P}$ will
be called a \emph{degree reducing interpolation basis} w.r.t. $\prec$ for (\ref{inter}).
Assume that $p_1 \prec p_2\prec \cdots \prec p_\mu$. If
$$
p_j(\xi^{(i)})=\delta_{ij},\quad 1\leq i \leq j \leq \mu,
$$
for some suitable reordering of $\Xi$, then we call $\{p_1, \ldots, p_\mu\}$ a
\emph{degree reducing interpolation Newton basis}(DRINB) w.r.t. $\prec$ for
(\ref{inter}).

Let $G_\prec$ be the reduced Gr\"{o}bner basis for the vanishing ideal $\mathcal{I}(\Xi)$
w.r.t.$\prec$. The set
$$
\mathrm{N}_\prec(\mathcal{I}(\Xi)):=\{x^{\bm{\alpha}}\in \mathbb{T}^2:\mathrm{LT}(g)\nmid
x^{\bm{\alpha}}, \forall g\in G_\prec\}
$$
is called the \emph{Gr\"{o}bner \'{e}scalier} of $\mathcal{I}(\Xi)$ w.r.t. $\prec$. From
\cite{Sau2006, dBo2007}, the interpolation space spanned by
$\mathrm{N}_\prec(\mathcal{I}(\Xi))$, denoted by $\mathcal{P}_\prec(\Xi)$, is canonical
since it is the unique degree reducing interpolation space spanned by monomials w.r.t.
$\prec$ for (\ref{inter}). Hence, we call $\mathrm{N}_\prec(\mathcal{I}(\Xi))$ the
\emph{degree reducing interpolation monomial basis}(DRIMB) w.r.t. $\prec$ for
(\ref{inter}), with $\#\mathrm{N}_\prec(\mathcal{I}(\Xi))=\mu$. Let
$$
\mathrm{N}_\prec(\Xi):=\{\bm{\alpha}: x^{\bm{\alpha}}\in
\mathrm{N}_\prec(\mathcal{I}(\Xi))\}\subset \mathbb{N}_0^2.
$$
We can deduce easily that $\mathrm{N}_\prec(\Xi)$ is a lower set and obviously has a
one-to-one correspondence with $\mathrm{N}_\prec(\mathcal{I}(\Xi))$. Therefore,
interpolation scheme $(\Xi, \mathcal{P}_\prec(\Xi))$ can be equivalently represented as
$(\Xi, \mathrm{N}_\prec(\Xi))$.

According to \cite{Cra2004}, we can construct two particular lower
sets from $\Xi$, denoted by $S_x(\Xi), S_y(\Xi)$, which reflect the
geometry of $\Xi$ in certain sense.

Specifically, we cover the points in $\Xi$ by lines $l_0^x, l_1^x,\ldots, l_\nu^x$
parallel to the $x$-axis and assume that, without loss of generality, there are $m_j+1$
points, say $u_{0j}^x, u_{1j}^x,\ldots, u_{m_j,j}^x$, on $l_j^x$ with $m_0\geq m_1\geq
\cdots\geq m_\nu\geq 0$ hence the ordinates of $u_{ij}^x$ and $u_{i'j}^x, i\neq i'$,
same. Now, we set
$$
  S_x(\Xi):=\{(i,j) : 0\leq i \leq m_j,\ 0\leq j\leq \nu\},
$$
which apparently equals to $\mathrm{L}_x(m_0,m_1,\ldots,m_\nu)$. We can also cover
 the points by lines $l_0^y, l_1^y,\ldots, l_\lambda^y$ parallel
to the $y$-axis and denote the points on line $l_i^y$ by $u_{i0}^y, u_{i1}^y,\ldots,
u_{i,n_i}^y$ with $n_0\geq n_1\geq \cdots\geq n_\lambda\geq 0$ hence the abscissae of
$u_{ij}^y$ and $u_{ij'}^y, j\neq j'$, same. Similarly, we put
$$
S_y(\Xi):=\{(i,j) : 0\leq i \leq \lambda,\ 0\leq j\leq n_i\}=\mathrm{L}_y(n_0, n_1,
\ldots, n_\lambda).
$$
In addition, we can also define the sets of abscissae and ordinates
\begin{equation}\label{HV}
    \begin{aligned}
H_j(\Xi):=&\{\bar{x}: (\bar{x}, \bar{y}) \in l_j^x\cap \Xi\},\quad 0\leq j \leq \nu,\\
V_i(\Xi):=&\{\bar{y}: (\bar{x}, \bar{y}) \in l_i^y\cap \Xi\},\quad 0\leq i \leq \lambda.
    \end{aligned}
\end{equation}

\begin{defn}\cite{Cra2004}
  We say that a set $\Xi$ of distinct points in $\mathbb{F}^2$ is \emph{cartesian} if there
  exists a lower set $\mathcal{A}$ such that $\Xi$ can be written as
$$
  \Xi=\{(x_i,y_j) : (i,j)\in \mathcal{A}\},
$$
  where the $x_i$'s are distinct numbers, and similarly the
  $y_j$'s. We also say that $\Xi$ is $\mathcal{A}$-cartesian.
\end{defn}

To the best of our knowledge, there are two criteria for determining whether a
2-dimensional point set is cartesian.

\begin{thm}\textup{\cite{Cra2004}}\label{Car}
A set of distinct points $\Xi\subset \mathbb{F}^2$ is cartesian if and only if
$S_x(\Xi)=S_y(\Xi)$.
\end{thm}

\begin{thm}\textup{\cite{CDZ2006}}\label{CDZ}
A set of distinct points $\Xi\subset \mathbb{F}^2$ is cartesian if and only if
$$
H_0(\Xi) \supseteq H_1(\Xi) \supseteq \cdots \supseteq H_\nu(\Xi),\quad
  V_0(\Xi) \supseteq V_1(\Xi) \supseteq \cdots \supseteq
  V_\lambda(\Xi).
$$

\end{thm}

About the bivariate Lagrange interpolation on a cartesian set, \cite{Cra2004} proved the
succeeding theorem.

\begin{thm}\textup{\cite{Cra2004}}\label{lowerset}
Given a cartesian set $\Xi\subset \mathbb{F}^2$, there exists a unique lower set
$\mathcal{A}\in \mathbb{N}_0^2$ such that $\Xi$ is $\mathcal{A}$-cartesian and the
Lagrange interpolation scheme $(\Xi, \mathcal{A})$ is regular.
\end{thm}

Finally, we will redescribe the classical BM algorithm with the notation established
above.

\begin{alg}\label{BMalg}
(BM Algorithm)

\textbf{Input}: A set of distinct points $\Xi=\{\xi^{(i)} : i=1, \ldots, \mu\}\subset
\mathbb{F}^d$ and
 a fixed term order $\prec$.\\
\indent \textbf{Output}: The 3-tuple $(G,N,Q)$, where $G$ is the reduced Gr\"{o}bner
basis for $\mathcal{I}(\Xi)$ w.r.t. $\prec$, $N$ is the Gr\"{o}bner \'{e}scalier of
$\mathcal{I}(\Xi)$ (the DRIMB for (\ref{inter}) also) w.r.t. $\prec$, and $Q$ is a DRINB
 w.r.t. $\prec$ for (\ref{inter}).
\\
\indent \textbf{BM1.} Start with lists $G=[\ ],N=[\ ],Q=[\ ], L=[1]$, and a matrix
$B=(b_{ij})$ over $\mathbb{F}$ with $\mu$ columns and zero rows initially.\\
\indent \textbf{BM2.} If $L=[\ ]$, return $(G,N,Q)$ and stop. Otherwise, choose the
monomial
$t=\mbox{min}_\prec L$, and delete $t$ from $L$.\\
\indent \textbf{BM3.} Compute the evaluation vector
$(t(\xi^{(1)}),\ldots,t(\xi^{(\mu)}))$, and reduce it against the rows of $B$ to obtain
$$
(v_1,\ldots,v_\mu)=(t(\xi^{(1)}),\ldots,t(\xi^{(\mu)}))-\sum_i
a_i(b_{i1},\ldots,b_{i\mu}), \quad a_i\in \mathbb{F}.
$$
\indent \textbf{BM4.}. If $(v_1,\ldots,v_\mu)=(0,\ldots,0)$, then append the polynomial
$t-\sum_i a_iq_i$ to the list $G$, where $q_i$ is the $i$th element of $Q$. Remove from
$L$ all the multiples of $t$.
Continue with \textbf{BM2}.\\
\indent \textbf{BM5.} Otherwise $(v_1,\ldots,v_\mu)\neq(0,\ldots,0)$, add
$(v_1,\ldots,v_\mu)$ as a new row to $B$ and $t-\sum_i a_iq_i$ as a new element to $Q$.
Append the monomial $t$ to $N$, and add to $L$ those elements of $\{x_1t,\ldots,x_dt\}$
that are neither multiples of an element of $L$ nor of $\mathrm{LT}(G)$. Continue with
\textbf{BM2}.
\end{alg}

\section{Special cases}\label{special}

In this section, we will focus on $\prec_{\mathrm{lex}}$ and $\prec_{\mathrm{inlex}}$
that may be the most talked about term orders. For these special cases, our preprocessing
paradigm will first provide exact $N, Q$ of the 3-tuple output $(G, N, Q)$ to BM
algorithm directly and effortlessly. And then, $G$ can be obtained by BM algorithm
easily. Note that we will continue with all the notation that we established for
$S_x(\Xi)$ and $S_y(\Xi)$ in the previous section.

\begin{prop}\label{phi^x}
Let $\Xi$ be a set of $\mu$ distinct points $u_{mn}^x=(x_{mn}, y_{mn})\in \mathbb{F}^2,
(m, n)\in S_x(\Xi)$. The points give rise to polynomials
\begin{equation}\label{phi}
\phi_{ij}^x=\varphi_{ij}^x \prod_{t=0}^{j-1}(y-y_{0t})\prod_{s=0}^{i-1}(x-x_{sj}), \quad
(i,j)\in S_x(\Xi),
\end{equation}
where
$\varphi_{ij}^x=1/\prod_{t=0}^{j-1}(y_{0j}-y_{0t})\prod_{s=0}^{i-1}(x_{ij}-x_{sj})\in
\mathbb{F}$, and the empty products are taken as 1. Then we have
$$
\phi_{ij}^x(u_{mn}^x)=\delta_{(i,j), (m,n)}, \quad (i,j)\succeq_{\mathrm{inlex}} (m,n).
$$
\end{prop}

\begin{pf}
Fix $(i, j)\in S_x(\Xi)$. Recalling the definition of $u_{ij}^x$, we have
$y_{0j}=y_{ij}$. If $(i, j)=(m, n)$, by $y_{00}\neq y_{01}\neq\cdots\neq y_{0j}$ and
$x_{0j}\neq x_{1j}\neq\cdots\neq x_{ij}$, we have
\begin{equation*}
\phi_{ij}^x(u_{ij}^x)=\varphi_{ij}^x\prod_{t=0}^{j-1}(y_{ij}-y_{0t})\prod_{s=0}^{i-1}(x_{ij}-x_{sj})
=\varphi_{ij}^x\prod_{t=0}^{j-1}(y_{0j}-y_{0t})\prod_{s=0}^{i-1}(x_{ij}-x_{sj}),
\end{equation*}
which implies $\phi_{ij}^x(u_{ij}^x)=1$.

Otherwise, if $(i, j)\succ_{\mathrm{inlex}} (m, n)$, we have $j>n$, or $j=n, i>m$. When
$j>n$, we have
\begin{align*}
\phi_{ij}^x(u_{mn}^x)&=\varphi_{ij}^x(y_{mn}-y_{00})\cdots(y_{mn}-y_{0n})\cdots(y_{mn}-y_{0,j-1})\prod_{s=0}^{i-1}(x_{mn}-x_{sj})\\
&=\varphi_{ij}^x(y_{0n}-y_{00})\cdots(y_{0n}-y_{0n})\cdots(y_{0n}-y_{0,j-1})\prod_{s=0}^{i-1}(x_{mn}-x_{sj})\\
&=0,
\end{align*}
and when $j=n, i>m$,
\begin{align*}
\phi_{ij}^x(u_{mn}^x)&=\varphi_{ij}^x\prod_{t=0}^{j-1}(y_{mn}-y_{0t})(x_{mn}-x_{0j})\cdots(x_{mn}-x_{mj})\cdots(x_{mn}-x_{i-1,j})\\
&=\varphi_{ij}^x\prod_{t=0}^{n-1}(y_{mn}-y_{0t})(x_{mn}-x_{0n})\cdots(x_{mn}-x_{mn})\cdots(x_{mn}-x_{i-1,n})\\
&=0,
\end{align*}
which leads to
$$
\phi_{ij}^x(u_{mn}^x)=0, \quad (i,j)\succ_{\mathrm{inlex}} (m,n).
$$
\qed

\end{pf}

Similarly, we can prove the following proposition:
\begin{prop}\label{phi^y}
Let $\Xi$ be a set of $\mu$ distinct points $u_{mn}^y=(x_{mn}, y_{mn})\in \mathbb{F}^2,
(m, n)\in S_y(\Xi)$. We define the polynomials
\begin{equation}\label{phiy}
\phi_{ij}^y=\varphi_{ij}^y\prod_{s=0}^{i-1}(x-x_{s0})\prod_{t=0}^{j-1}(y-y_{it}),\quad
(i,j)\in S_y(\Xi),
\end{equation}
where
$\varphi_{ij}^y=1/\prod_{s=0}^{i-1}(x_{i0}-x_{s0})\prod_{t=0}^{j-1}(y_{ij}-y_{it})\in
\mathbb{F}$. The empty products are taken as 1. Then,
$$
\phi_{ij}^y(u_{mn}^y)=\delta_{(i,j),(m,n)}, \quad (i,j)\succeq_{\mathrm{lex}} (m,n).
$$
\end{prop}

In 2004, \cite{Cra2004} proved that the Lagrange interpolation schemes $(\Xi, S_x(\Xi))$
and $(\Xi, S_y(\Xi))$ are both regular. Here we reprove the regularities in another way
for the purpose of presenting the degree reducing interpolation bases theoretically .

\begin{thm}\label{minbasis}
Resume the notation in \emph{Proposition \ref{phi^x}} and
\emph{\ref{phi^y}}. Then the Lagrange interpolation schemes $(\Xi,
S_x(\Xi))$ and $(\Xi, S_y(\Xi))$ are
regular. Furthermore,\\
\textup{(i)} the set $N_x:=\{x^iy^j : (i,j)\in S_x(\Xi)\}$ is the
\emph{DRIMB} as well as $Q_x:=\{\phi_{ij}^x: (i, j)\in S_x(\Xi)\}$
is a \emph{DRINB}  w.r.t. $\prec_{\mathrm{lex}}$
for the interpolation problem \textup{(\ref{inter})}.\\
\textup{(ii)} the set $N_y:=\{x^iy^j : (i, j)\in S_y(\Xi)\}$ is the \emph{DRIMB} as well
as $Q_y:=\{\phi_{ij}^y: (i, j)\in S_y(\Xi)\}$ is a \emph{DRINB} w.r.t.
$\prec_{\mathrm{inlex}}$ for \textup{(\ref{inter})}.
\end{thm}

\begin{pf}
We only give the proof for $S_x(\Xi)$. The statements about
$S_y(\Xi)$ can be proved likewise.

First, we will show the regularity of the interpolation scheme
$(\Xi, S_x(\Xi))$. Let $
\mathcal{P}_x:=\mathrm{Span}_{\mathbb{F}}N_x\subset \Pi^2$ with
$\dim\mathcal{P}_x=\#\Xi=\mu$. Obviously, $N_x$ is the monomial
basis for it. By (\ref{phi}), we can check easily that
$$
\mathrm{Span}_\mathbb{F}Q_x\subseteq \mathcal{P}_x.
$$
Construct a square matrix $B_{\mu\times\mu}$ whose $(h, k)$ entry is $\phi_h^x(u_k^x)$
where $\phi_h^x, u_k^x$ are $h$th and $k$th elements of $Q_x$ and $\Xi=\{u_{mn}^x: (m,
n)\in S_x(\Xi)\}$ w.r.t. the increasing $\prec_{\mathrm{inlex}}$ on $(i, j)$ and $(m, n)$
respectively. From Proposition \ref{phi^x}, $B_{\mu\times\mu}$ is upper unitriangular
which implies that $\mathrm{Span}_\mathbb{F}Q_x=\mathcal{P}_x$ and $Q_x$ forms a Newton
basis for $\mathcal{P}_x$. It follows that $\mathcal{P}_x$ is an interpolation space for
Lagrange interpolation (\ref{inter}) therefore the scheme $(\Xi, \mathcal{P}_x)$ is
regular. Since $(\Xi, S_x(\Xi))=(\Xi, \mathcal{P}_x)$, according to Section \ref{pre},
$(\Xi, S_x(\Xi))$ is regular.

Next, we shall verify that the statements in (i), which is equivalent to the statement
that $\mathcal{P}_x$ is a degree reducing interpolation space w.r.t.
$\prec_{\mathrm{lex}}$ for (\ref{inter}) that coincides with
$\mathcal{P}_{\prec_{\mathrm{lex}}}(\Xi)$. Since the arguments above have proved that
$\mathcal{P}_x$ satisfies the \textbf{DR1} condition in Definition \ref{MDIS}, what is
left for us is to check the \textbf{DR2} condition. From \cite{dBo2007}, we only need to
check it for monomials.

Take a monomial $x^{i_0}y^{j_0}\in\mathbb{T}^2$. We shall prove that
\begin{equation}\label{DR2}
  L_{\mathcal{P}_x}x^{i_0}y^{j_0} \preceq_{\mathrm{lex}} x^{i_0}y^{j_0}.
\end{equation}
Since $\mathcal{P}_x$ satisfies \textbf{DR1}, $L_{\mathcal{P}_x}x^{i_0}y^{j_0}$ is the
unique polynomial in $\mathcal{P}_x$ that matches $x^{i_0}y^{j_0}$ on $\Xi$. Therefore,
when $x^{i_0}y^{j_0}\in N_x$, we have $L_{\mathcal{P}_x}x^{i_0}y^{j_0}=x^{i_0}y^{j_0}$ ,
namely (\ref{DR2}) is true for this case. Assume that
$$
S_x(\Xi)=\mathrm{L}_x(m_0, \ldots, m_{n_0})=\mathrm{L}_y(n_0,
\ldots, n_{m_0}).
$$
It is easy to see that $\delta(\mathcal{P}_x)=(m_0, n_{m_0})$. If
$x^{m_0}y^{n_{m_0}}\prec_{\mathrm{lex}}x^{i_0}y^{j_0} $ then
$\delta(L_{\mathcal{P}_x}x^{i_0}y^{j_0})\preceq_{\mathrm{lex}}
\delta(\mathcal{P}_x)=(m_0, n_{m_0})\prec_{\mathrm{lex}}(i_0,
j_0)=\delta(x^{i_0}y^{j_0})$ that leads to (\ref{DR2}) for the case.

Thus, what remains for us is to check (\ref{DR2}) for $x^{i_0}y^{j_0}\notin N_x$ with
$(i_0, j_0)\prec_{\mathrm{lex}}(m_0, n_{m_0})$ that implies $0\leq i_0<m_0, j_0>n_{i_0}$.
 For this, we only need to verify that

\begin{equation}\label{Lp}
L_{\mathcal{P}_x}x^{i_0}y^{j_0}\in \mathrm{Span}_\mathbb{F}\{x^iy^j:
(i, j)\in F_{i_0}\},
\end{equation}
where $ F_{i_0}=\{(i, j)\in S_x(\Xi): (i,j)\prec_{\mathrm{lex}}(i_0,
j_0)\}\subset S_x(\Xi). $ If $x^{i_0}y^{j_0}\in \mathcal{I}(\Xi)$,
then $L_{\mathcal{P}_x}x^{i_0}y^{j_0}=0\prec_{\mathrm{lex}}
x^{i_0}y^{j_0}$. The statement (\ref{Lp}) becomes trivial in this
case. Otherwise, if we can find a polynomial $p\in \Pi^2$ such that
\begin{eqnarray}\label{pol}
p=x^{i_0}y^{j_0}-\sum_{(i,j)\in F_{i_0}}a_{ij}x^iy^j \in
\mathcal{I}(\Xi),
\end{eqnarray}
where $a_{ij}\in \mathbb{F}$ are not all zero, then (\ref{Lp}) follows.

According to Section \ref{pre}, our point set
$\Xi=\{u_{ij}^x=(x_{ij}, y_{ij}): (i, j)\in S_x(\Xi)\}$. Let
$\Xi'=\{u_{mn}^x\in \Xi: (m, n)\in F_{i_0}\}\subset \Xi$. Now, we
claim that there exists a unique polynomial $p$ of the form
(\ref{pol}) such that $p\in \mathcal {I}(\Xi')$, which is equivalent
to the statement that the linear system
\begin{eqnarray}\label{ls}
\sum_{(i,j)\in F_{i_0}}a_{ij}x_{mn}^iy_{mn}^j= x_{mn}^{i_0}y_{mn}^{j_0},\quad u_{mn}^x\in
\Xi',
\end{eqnarray}
has a unique solution.

Note that $\mathrm{Span}_{\mathbb{F}}\{ x^iy^j : (i,j)\in
F_{i_0}\}=\mathrm{Span}_{\mathbb{F}}\{ \phi^x_{ij} : (i,j)\in F_{i_0}\}$. We can conclude
that the rank of the coefficient matrix of (\ref{ls}) is equal of that of the matrix
$B'_{\#F_{i_0}\times\#F_{i_0}}$, which is a submatrix of $B$ whose $(h, k)$ entry is
$\phi_h^x(u_k^x)$ where $\phi_h^x, u_k^x$ are $h$th and $k$th elements of $\{\phi_{ij}^x:
(i, j)\in F_{i_0}\}$ and $\Xi'=\{u_{mn}^x\}$ w.r.t. the increasing
$\prec_{\mathrm{inlex}}$ on $(i, j)$ and $(m, n)$ respectively. By (\ref{phi}), we see
easily that $B'$ is upper unitriangular which implies that the coefficient matrix of
(\ref{ls}) is of full rank. Accordingly, there is a unique polynomial $p\in \mathcal
{I}(\Xi')$ that has the form (\ref{pol}).

Now we shall verify that $p(u_{ij}^x)=0, u_{ij}^x\in\Xi\setminus\Xi'$. By the definition
of $\Xi'$, we know that $i>i_0$ here. Let
$$
q(x):=p(x,y_{ij})=\sum_{s=0}^{i_0}b_sx^s\in \Pi^1,\quad b_s\in\mathbb{F}.
$$
Since  $y_{0j}=y_{1j}=\cdots=y_{i_0,j}=y_{ij}$ and
$u_{0j}^x,u_{1j}^x,\ldots,u_{i_0,j}^x\in \Xi'$, it follows that
$$
q(x_{sj})=p(x_{sj}, y_{ij})=p(x_{sj}, y_{sj})=p(u_{sj}^x)=0,\quad s=0,\ldots,i_0,
$$
namely $q(x)$ has $i_0+1$ zero points that clearly implies $q(x)\equiv 0$. Since
$p(u_{ij}^x)=q(x_{ij})=0$, we have $p\in \mathcal{I}(\Xi)$. By (\ref{Lp}), (\ref{DR2}) is
true in this case. As a result, for any $f\in \Pi^2$, we have
$$
L_{\mathcal{P}_x}f \preceq_{\mathrm{lex}} f,
$$
that is to say $\mathcal{P}_x$ satisfies \textbf{DR2}.

Consequently, by Definition \ref{MDIS}, $\mathcal{P}_x$ is a degree reducing
interpolation space w.r.t. $\prec_{\mathrm{lex}}$ for Lagrange interpolation
(\ref{inter}). Hence $N_x$ is the DRIMB and $Q_x$ is a Newton basis w.r.t.
$\prec_{\mathrm{lex}}$ for (\ref{inter}).\qed

\end{pf}

Note that  $\mathcal{P}_{\prec_{\mathrm{lex}}}(\Xi)$ is the unique degree reducing
interpolation space spanned by monomials w.r.t. $\prec_{\mathrm{lex}}$, thus we have
$\mathcal{P}_x=\mathcal{P}_{\prec_{\mathrm{lex}}}(\Xi)$. Therefore,
$N_x=N_{\prec_{\mathrm{lex}}}(\mathcal {I}(\Xi))$ holds, which means that $N_x$ is also
the Gr\"{o}bner  \'{e}scalier of $\mathcal {I}(\Xi)$ w.r.t. $\prec_{\mathrm{lex}}$.

\begin{cor}\label{3lower}
  If $\Xi\subset \mathbb{F}^2$ is an $\mathcal{A}$-cartesian set, then $\mathcal{A}=S_x(\Xi)=S_y(\Xi)$.
\end{cor}

\begin{pf}
  Since $\Xi$ is cartesian, by Theorem \ref{Car} and \ref{minbasis}, we have
  $S_x(\Xi)=S_y(\Xi)$ hence $(\Xi, S_x(\Xi))=(\Xi, S_y(\Xi))$ are both regular. But from Theorem \ref{lowerset},
  only $\mathcal{A}$ can make $(\Xi, \mathcal{A})$ regular, therefore $\mathcal{A}=S_x(\Xi)=S_y(\Xi)$.\qed
\end{pf}

From Algorithm \ref{BMalg} we know that $G, N, Q$ are essential elements of BM algorithm
and compose its output. For $\prec_{\mathrm{lex}}$ and $\prec_{\mathrm{inlex}}$ cases,
Theorem \ref{minbasis} presents us $N$ and $Q$ theoretically hence we can obtain them
with little effort. According to \cite{MMM1993}, the leading terms of $G$ are contained
in the border set of $N$. Therefore, we can get $G$ faster than compute $G$ directly with
BM algorithm. Now is our algorithm.

\begin{alg}\label{pBMlex}(SPBM)

\textbf{Input}: A set of distinct affine points $\Xi\subset \mathbb{F}^2$ and fixed
$\prec_{\mathrm{lex}}$ or $\prec_{\mathrm{inlex}}$.

 \textbf{Output}: The 3-tuple $(G,N,Q)$,
where $G$ is the reduced Gr\"{o}bner basis of $\mathcal{I}(\Xi)$, $N$ is the Gr\"{o}bner
\'{e}scalier $\mathrm{N}(\mathcal{I}(\Xi))$, and $Q$
is a DRINB for the Lagrange interpolation on $\Xi$.\\
\indent \textbf{SPBM1.} Construct lower set $S_x(\Xi)$ or $S_y(\Xi)$
according to Section \ref{pre}.
\\
\indent \textbf{SPBM2.} Compute the sets $N$ and $Q$ by Theorem \ref{minbasis}.
\\
\indent \textbf{SPBM3.} Construct the border set $L:=\{x\cdot t : t\in N\}\bigcup\{
y\cdot t : t\in N\}\setminus N$ and the matrix $B$ that is same to the $B_{\mu\times\mu}$
in the proof of Theorem \ref{minbasis}.
\\
\indent \textbf{SPBM4.} Goto \textbf{BM2} of BM algorithm for the reduced Gr\"{o}bner
basis $G$.
\end{alg}

\begin{exmp}\label{eg21}
Let
$$
\Xi=\{(0,1),(0,3),(1,0),(1,2),(1,3),(1,4),(2,1),(2,2),(3,1)\}\subset \mathbb{Q}^2.
$$
First, we choose lines $x=1, x=0, x=2, x=3$ as $l_0^y, l_1^y, l_2^y, l_3^y$ respectively
\textup{(}Shown in (a) of \textup{Figure \ref{eg2})}, therefore we have
$$
S_y=\{(0,0),(0,1),(0,2),(0,3),(1,0),(1,1),(2,0),(2,1),(3,0)\},
$$
which is illutrated in (b) of \textup{Figure \ref{eg2}}.

\begin{figure}[!htbp]
\centering
\subfigure[$\Xi$]{\includegraphics[width=6cm,height=6cm]{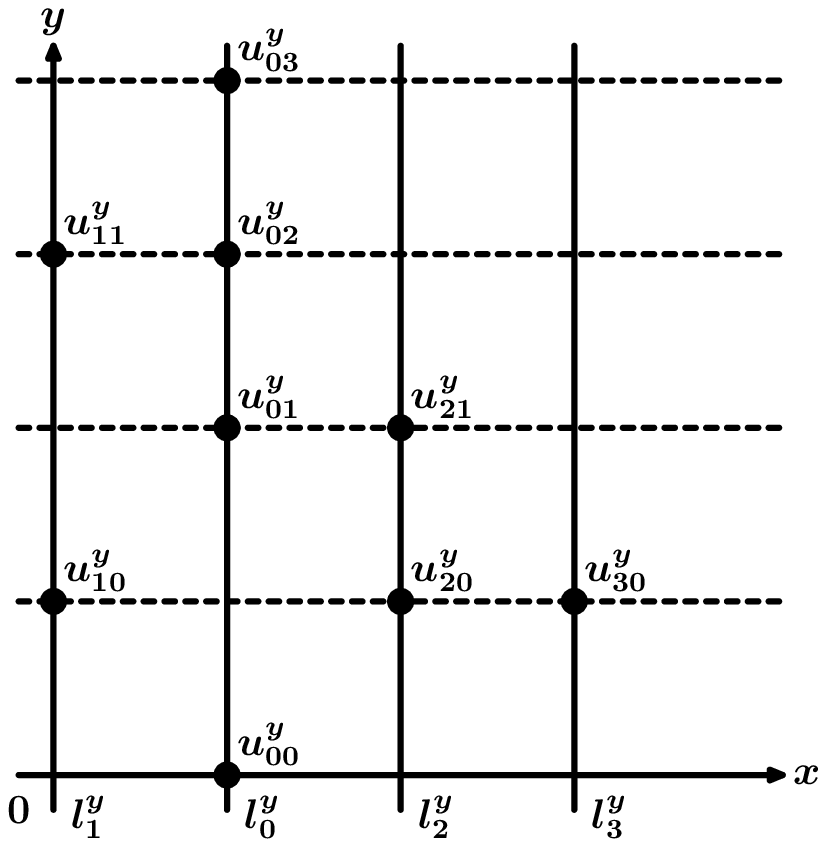}}
\subfigure[$S_y$]{\includegraphics[width=6cm,height=6cm]{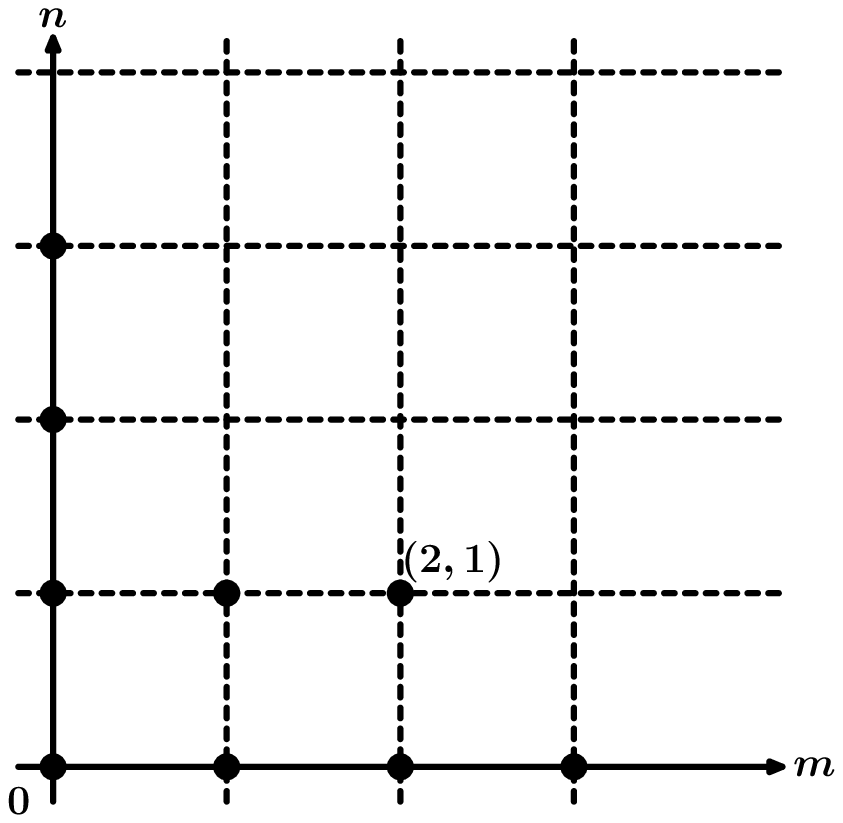}}
\caption{The point
set and related $S_y$ of Example \ref{eg21}.}\label{eg2}
\end{figure}
Thus, by Theorem \ref{minbasis}, we have
\begin{align*}
N=\{&1,y,y^2,y^3,x,xy,x^2,x^2y,x^3\}; \\
Q=\{&1, \frac{1}{2}y, \frac{1}{3}y^2-\frac{2}{3}y, \frac{1}{8}y^3-\frac{5}{8}y^2+\frac{3}{4}y,
    -x+1, -\frac{1}{2}xy+\frac{1}{2}y+\frac{1}{2}x-\frac{1}{2},\\
    & \frac{1}{2}x^2-\frac{1}{2}x, \frac{1}{2}x^2y-\frac{1}{2}xy-\frac{1}{2}x^2+\frac{1}{2}x,
    \frac{1}{6}x^3-\frac{1}{2}x^2+\frac{1}{3}x\}.\\
\end{align*}
Next, from \textbf{SPBM3}, the border set $ L=\{y^4, xy^2, xy^3, x^2y^2, x^3y, x^4\}$ and
the matrix
$$
B=\left(
    \begin{array}{cccc}
      1 & 1 & 1 & \cdots \\
      0 & 1 & 3/2 & \cdots \\
      0 & 0 & 1 & \cdots \\
      \vdots & \vdots & \vdots & \ddots \\
    \end{array}
  \right).
$$
Finally, turn to \textbf{BM2} with these $N, Q, L, B$, we can get the reduced Gr\"{o}bner
basis
\begin{align*}
G=\{&x^4-6x^3+11x^2-6x, x^3y-3x^2y+2xy-x^3+3x^2-2x,\\
&xy^2-y^2+\frac{1}{2}x^2y-\frac{9}{2}xy+4y-\frac{1}{2}x^2+\frac{7}{2}x-3,\\
&y^4-9y^3+26y^2-\frac{9}{2}x^2y+\frac{15}{2}xy-27y-3x^3+\frac{39}{2}x^2-\frac{51}{2}x+9\}.
\end{align*}
for $\mathcal{I}(\Xi)$ w.r.t. $\prec_{\mathrm{inlex}}$.
\end{exmp}

\begin{exmp}\label{eg11}
Given a bivariate point set
$$
\Xi=\{(0, 0), (0, 2), (0, 3), (1, 1), (\frac{5}{2}, 0), (\frac{5}{2}, 1), (\frac{5}{2},
2), (4, 0), (4, 2)\}\subset \mathbb{Q}^2.
$$
We choose lines $y=0, y=2, y=1, y=3$ as $l_0^x, l_1^x, l_2^x, l_3^x$ respectively
(Illustrated in (a) of Figure \ref{eg1}),  which follows that
$$
S_x=\{(0,0),(1,0),(2,0),(0,1),(1,1),(2,1),(0,2),(1,2),(0,3)\}.
$$

\begin{figure}[!htbp]
\centering \subfigure[$\Xi$]{\includegraphics[width=6cm,height=6cm]{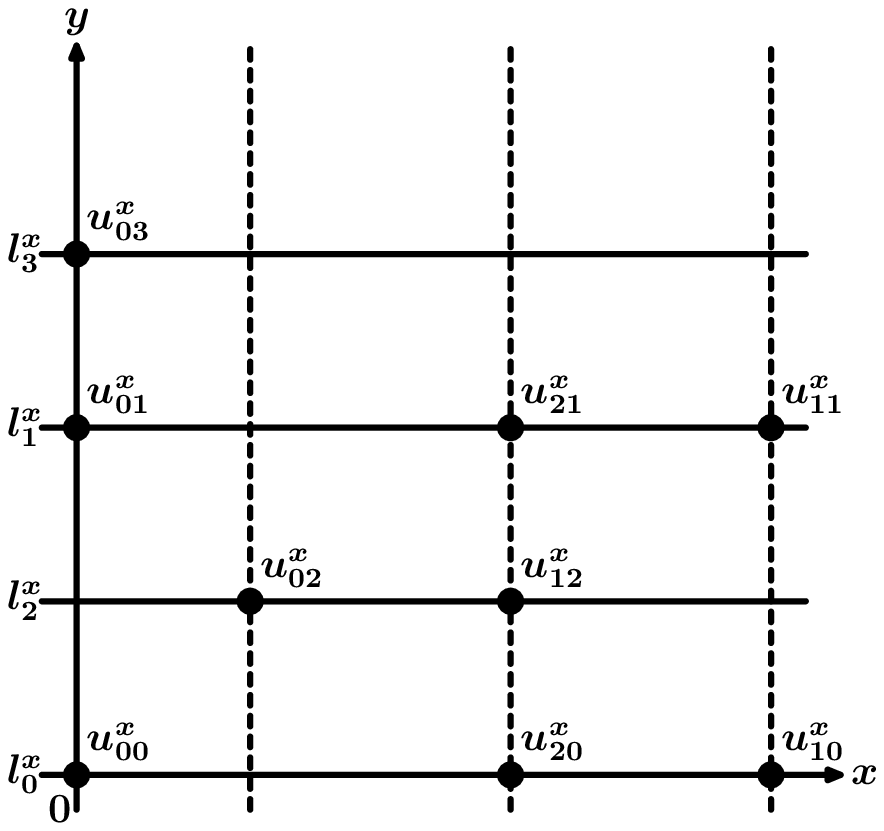}}
\subfigure[$S_x$]{\includegraphics[width=6cm,height=6cm]{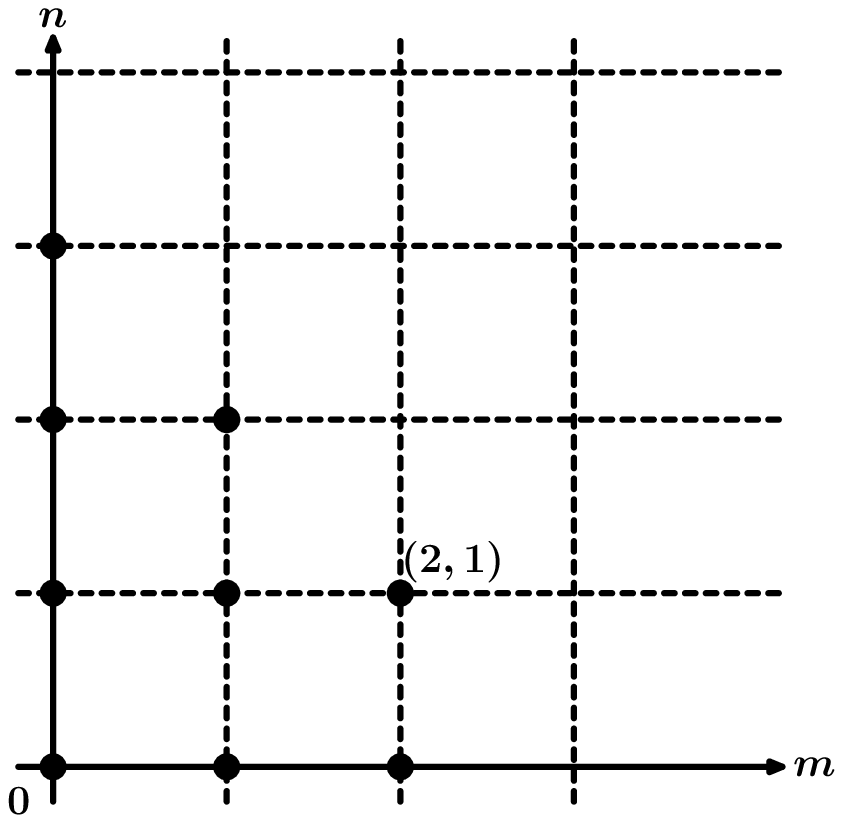}}
\caption{Illustrations for Example \ref{eg11}.}\label{eg1}
\end{figure}
Thus, with SPBM algorithm, we have
\begin{align*}
N=\{&1,x,x^2,y,xy,x^2y,y^2,xy^2,y^3\}, \\
Q=\{&1, \frac{1}{4}x, -\frac{4}{15}x^2+\frac{16}{15}x, \frac{1}{2}y, \frac{1}{8}xy, -\frac{2}{15}x^2y+\frac{8}{15}xy, -y^2+2y,\\
    & -\frac{2}{3}xy^2+\frac{2}{3}y^2+\frac{4}{3}xy-\frac{4}{3}y,
    \frac{1}{6}y^3-\frac{1}{2}y^2+\frac{1}{3}y\},\\
G=\{&y^4-6y^3+11y^2-6y, xy^3-3xy^2+2xy,
x^2y^2-2x^2y-\frac{7}{2}xy^2\\
&+7xy-\frac{5}{4}y^3+\frac{25}{4}y^2-\frac{15}{2}y, x^3-\frac{13}{2}x^2-3xy^2+6xy\\
&+10x-\frac{15}{4}y^3+\frac{75}{4}y^2-\frac{45}{2}y\}.
\end{align*}
\end{exmp}

\section{General cases}\label{general}

Next, we will discuss how to accelerate BM algorithm with respect to term orders other
than $\prec_{\mathrm{lex}}$ or $\prec_{\mathrm{inlex}}$.
 In \cite{Cra2004}, the author proposed that if the set of
points $\Xi$ is cartesian, then we can obtain the interpolation basis without any
difficulty, see Theorem \ref{lowerset}. But in general $\Xi$ may not be cartesian.
However, we have the following proposition.

\begin{prop}\label{lowersub}
There must exist at least one cartesian subset for any non-empty set of points in
$\mathbb{F}^2$.
\end{prop}

\begin{pf}
Let $\Xi$ be a non-empty set of points. Hence, there exists at least one point $\xi\in
\Xi$. But $\xi$ itself can construct a cartesian subset $\{\xi\}\subset \Xi$.\qed
\end{pf}

\begin{defn}
Let $\Xi$ be a set of points in $\mathbb{F}^2$ and $\Xi'$ be a cartesian subset of $\Xi$.
We say that $\Xi'$ is a \emph{maximal cartesian subset} of $\Xi$ if any cartesian proper
subset $\Xi''$ of $\Xi$ containing $\Xi'$ is such that $\Xi''=\Xi'$. In addition, a
\emph{maximal row subset} of $\Xi$ is a non-empty subset that equals the intersection of
$\Xi$ and a horizontal line.
\end{defn}

From Proposition \ref{lowersub} we know that, for a set of given points, we can surely
find a maximal cartesian subset of it. Is it unique? Unfortunately, the answer is often
false.

\begin{exmp}\label{eg111}
Recall Example \ref{eg11}, let
\begin{eqnarray*}
\Xi'_1&=&\{(0,0),(0,2),(\frac{5}{2},0),(\frac{5}{2},1),(\frac{5}{2},2),(4,0),
(4,2)\},\\
\Xi'_2&=&\{(0,0),(0,2),(0,3),(\frac{5}{2},0),(\frac{5}{2},2),(4,0),
(4,2)\},\\
\Xi'_3&=&\{(1,1),(\frac{5}{2},0),(\frac{5}{2},1),(\frac{5}{2},2)\}.
\end{eqnarray*}
We can check easily that $\Xi_1',\Xi_2',\Xi_3'$ are all maximal cartesian subsets of
$\Xi$ (Illustrated in Figure \ref{eg1S}).

\begin{figure}[!htp]
\centering
\subfigure[$\Xi'_1$]{\includegraphics[width=6cm,height=6cm]{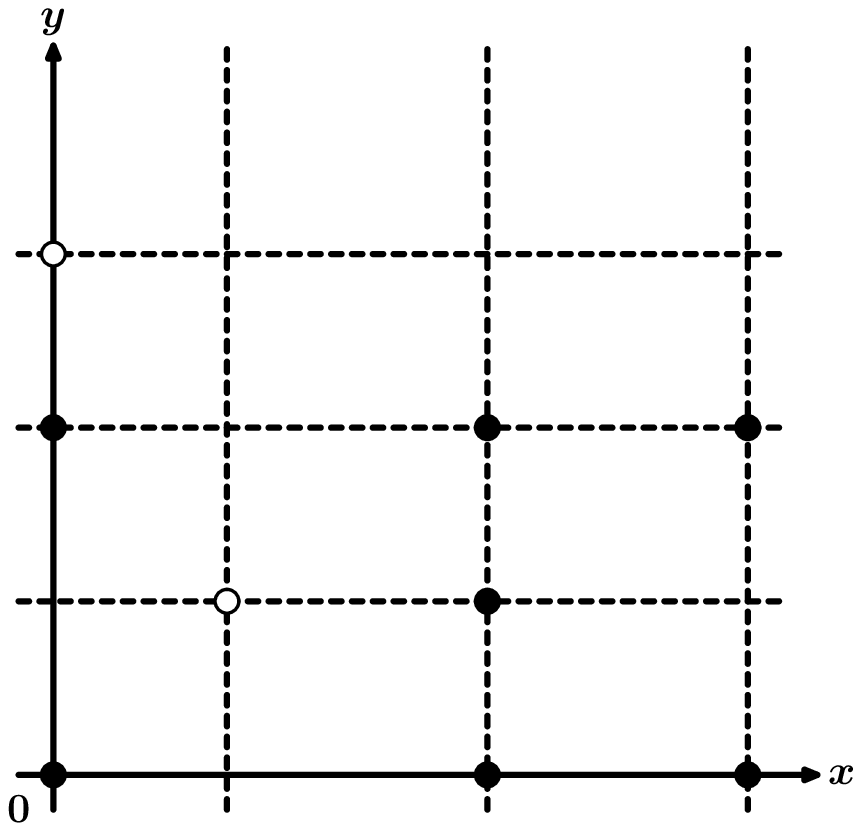}}
\subfigure[$\Xi'_2$]{\includegraphics[width=6cm,height=6cm]{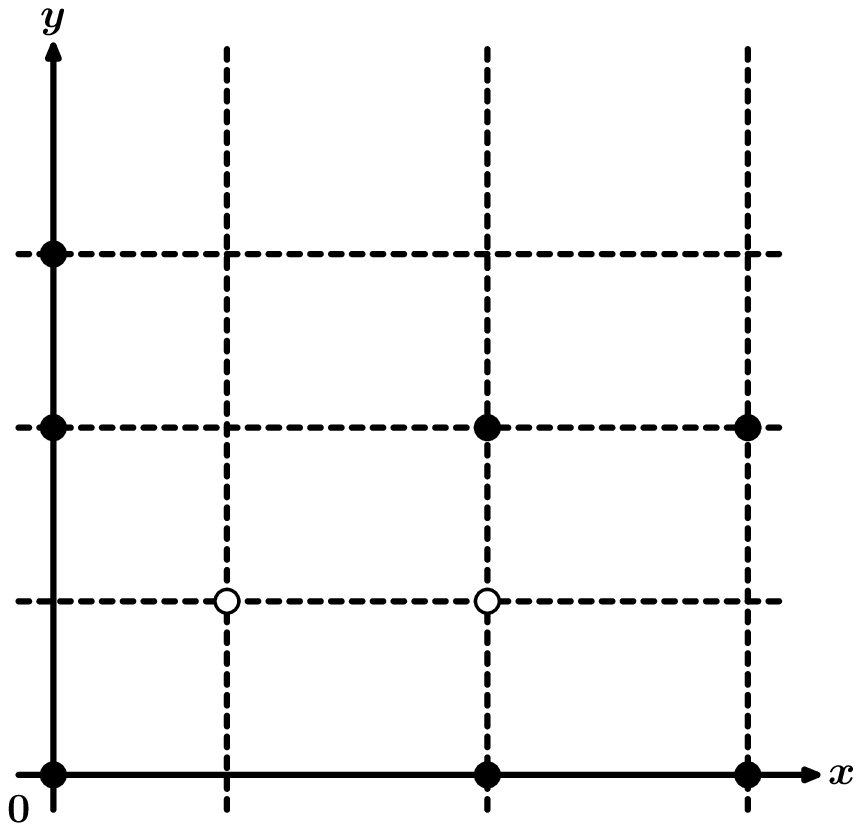}}\\
\subfigure[$\Xi'_3$]{\includegraphics[width=6cm,height=6cm]{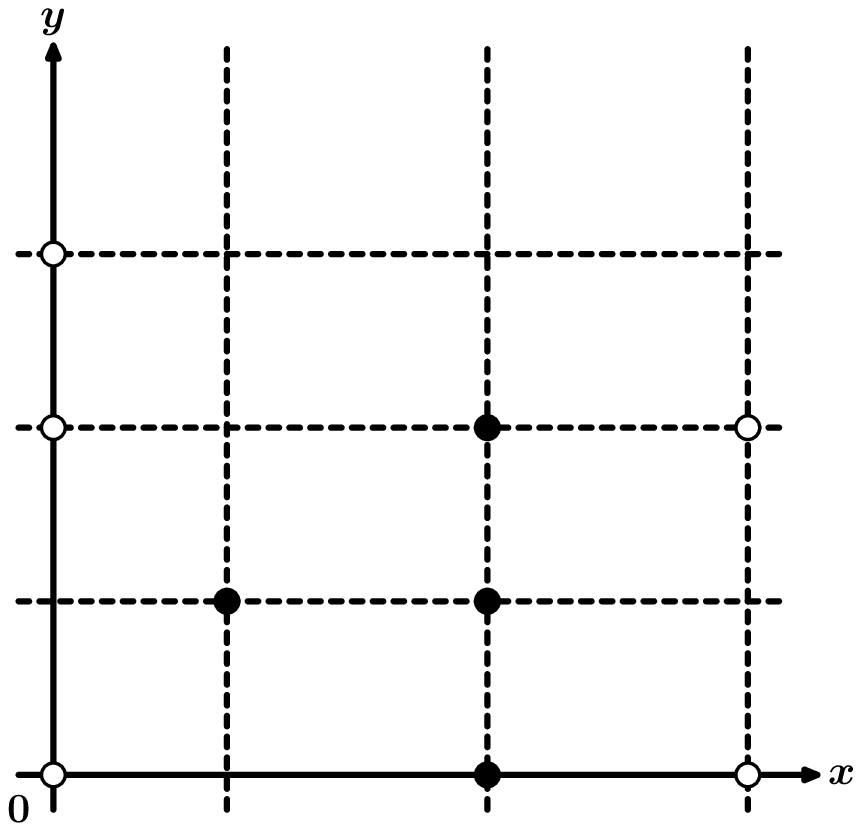}} \caption{Maximal
cartesian subsets of $\Xi$, where $\bullet$ denotes the points in $\Xi'_i, i=1,2,3$,
while $\circ$ denotes the points in $\Xi\backslash\Xi'_i$.}\label{eg1S}
\end{figure}
\end{exmp}

\begin{lem}\label{subcon}
Let $\Xi$ be a set of distinct points in $\mathbb{F}^2$ and $\prec$ a fixed term order.
If $\Xi'$ is an $\mathcal{A}'$-cartesian subset of $\Xi$, then
$$
\mathcal{A}' = \mathrm{N}_\prec(\Xi') \subset \mathrm{N}_\prec(\Xi),
$$
or equivalently,
$$
\{x^iy^j: (i, j)\in \mathcal{A}'\} = \mathrm{N}_\prec(\mathcal{I}(\Xi')) \subset
\mathrm{N}_\prec(\mathcal{I}(\Xi)).
$$
\end{lem}

\begin{pf}
From Section \ref{pre}, the Gr\"{o}bner \'{e}scalier
$\mathrm{N}_\prec(\mathcal{I}(\Xi'))$ is the DRIMB w.r.t. $\prec$ for the bivariate
Lagrange interpolation on $\Xi'$ hence the interpolation scheme $(\Xi',
\mathrm{N}_\prec(\Xi'))$ is regular. Since $\mathcal{A}'\subset \mathbb{N}_0^2$ is lower
and $\Xi'$ is $\mathcal{A}'$-cartesian, according to Theorem \ref{lowerset},
$\mathcal{A}'$ is the unique lower set making the bivariate Lagrange interpolation on
$\Xi'$ regular. This gives
$$
\mathcal{A'} = \mathrm{N}_\prec(\Xi').
$$
Since $\Xi'\subset\Xi$, from \cite{BW1993}, we know that the vanishing ideals satisfy
$\mathcal{I}(\Xi')\supset \mathcal{I}(\Xi)$. Denote by $G', G$ the reduced Gr\"{o}bner
bases for $\mathcal{I}(\Xi')$ and $\mathcal{I}(\Xi)$ w.r.t. $\prec$ respectively. We will
prove $\mathrm{N}_\prec(\mathcal{I}(\Xi')) \subset \mathrm{N}_\prec(\mathcal{I}(\Xi))$ by
contradiction. For any $x^iy^j\in \mathrm{N}_\prec(\mathcal{I}(\Xi'))$, we suppose there
were some $g\in G$ such that $\mathrm{LT}(g)| x^i y^j$. By \cite{BW1993},
$$
\langle \mathrm{LT}(G') \rangle = \langle \mathrm{LT}(\mathcal{I}(\Xi') \rangle \supset
\mathrm{LT}(\mathcal{I}(\Xi))\supset \mathrm{LT}(G).
$$
Therefore, $\mathrm{LT}(g)\in \mathrm{LT}(G) \subset \langle \mathrm{LT}(G') \rangle$
implies that there exists some $g'\in G'$ such that $\mathrm{LT}(g')| \mathrm{LT}(g)$.
Since $\mathrm{LT}(g)| x^i y^j$, we have $\mathrm{LT}(g')| x^i y^j$ that contradicts our
assumption on $x^iy^j$, which proves that $\mathrm{N}_\prec(\mathcal{I}(\Xi')) \subset
\mathrm{N}_\prec(\mathcal{I}(\Xi))$ due to the definition of
$\mathrm{N}_\prec(\mathcal{I}(\Xi))$. Finally, $\mathrm{N}_\prec(\Xi')\cong
\mathrm{N}_\prec(\mathcal{I}(\Xi'))$ and $\mathrm{N}_\prec(\Xi)\cong
\mathrm{N}_\prec(\mathcal{I}(\Xi))$ complete the proof. \qed

\end{pf}

\begin{rmk}
  For any $\mathcal{A}$-cartesian set $\Xi$, by Corollary \ref{3lower}, we have
$\mathcal{A}=S_x(\Xi)=S_y(\Xi)$ that obviously leads to
$\mathcal{A}=S_x(\Xi)=S_y(\Xi)=\mathrm{N}_\prec(\Xi)$, according to the Lemma above,
where term order $\prec$ is arbitrary.
\end{rmk}

Now comes an algorithm for constructing a maximal cartesian subset of a given point set
in $\mathbb{F}^2$.

\begin{alg}\label{MaxSubset}(Maximal Cartesian Subset Construction Algorithm)
\vskip 3mm
\textbf{Input}: A set of distinct points $\Xi=\{\xi^{(i)} : i=1, \ldots, \mu\}\subset \mathbb{F}^2$.\\
\indent\textbf{Output}: A maximal cartesian subset $\Xi'$ of $\Xi$.
\\
\indent\textbf{MCS1}. Start with an empty list $\Xi'=[\ ]$.\\
\indent\textbf{MCS2}. If $\Xi=[\ ]$, return the set $\Xi'$
and stop. Otherwise, compute lower sets $S_x(\Xi)$ and $S_y(\Xi)$.\\
\indent\textbf{MCS3}. If $S_x(\Xi)=S_y(\Xi)$, then replace $\Xi'$ by $\Xi'\cup \Xi$,
return the set $\Xi'$ and stop.
\\
\indent\textbf{MCS4}. Otherwise, we first choose a maximal row subset of $\Xi$ with
maximal cardinal number, denoted by $A$. Next,  delete from $\Xi$ the points either in
$A$ or have different abscissae from the points in $A$. Finally, replace $\Xi'$ by
$\Xi'\cup A$ and continue with \textbf{MCS2}.
\end{alg}

The following theorem ensure that this algorithm will terminate in finite steps with a
maximal cartesian subset as its output.

\begin{thm}
The algorithm described above will stop in a finite number of loops. Furthermore, the set
$\Xi'$ returned by the algorithm is a maximal cartesian subset.
\end{thm}

\begin{pf}
As input data of Algorithm \ref{MaxSubset}, point set $\Xi$ is finite. Observing that
$\#\Xi$ decreases actually in every loop, the algorithm will terminate in a finite
number, say $M$, of loops for sure. We assume that $M>1$ since $M=1$ is trivial.

$\Xi'_{\mathrm{in}}$ and $\Xi'_{\mathrm{out}}$ signify the input and output $\Xi'$ of
\textbf{MCS4} in some loop respectively. Next, we will prove by induction on $1\leq r
\leq M-1$ that in the $r$th loop $\Xi'_{\mathrm{out}}$ is a cartesian set. The case $r=1$
is obvious since $\Xi'_{\mathrm{in}}=[\ ]$ and $\Xi'_{\mathrm{out}}$ is clearly cartesian
as a maximal row subset of $\Xi$. Assume the statement is true for $r=l<M-1$. When
$r=l+1$, by the induction hypothesis, $\Xi'_{\mathrm{in}}$ is cartesian. Therefore, by
Corollary \ref{3lower}, we assume that
$$
\Xi'_{\mathrm{in}}=\{(x_i,y_j) : (i,j)\in S_x(\Xi'_{\mathrm{in}})\},
$$
where $S_x(\Xi'_{\mathrm{in}})=\mathrm{L}_x(m_0, \ldots, m_{n_0})=\mathrm{L}_y(n_0,
\ldots, n_{m_0})$. Observing the construction process of $\Xi'$ in the algorithm, we see
easily that $n_0=n_1=\cdots=n_{m_{n_0}}$. Let the maximal row subset of $\Xi$ we choose
at this moment be
$A=\{(\overline{x}^{(0)},\overline{y}),(\overline{x}^{(1)},\overline{y}),\ldots,(\overline{x}^{(k)},\overline{y})\}$.
Due to the nature of $A$,  we have $k\leq m_{n_0}$ and $\overline{y}\neq y_j, j=0,\ldots,
n_0$.

We claim that the set $\Xi'_{\mathrm{in}}\cup A$ is cartesian. In fact, we will focus on
the horizontal parallel lines $l_j^x: y=y_j, j=0,\ldots,n_0,$ and $l_{n_0+1}^x:
y=\overline{y}$. Resume the notation in (\ref{HV}). $H_j(\Xi'_{\mathrm{in}}\cup
A)=H_j(\Xi'_{\mathrm{in}})=\{x_i : 0\leq i\leq m_j\}, j=0,\ldots, n_0$, and
$H_{n_0+1}(\Xi'_{\mathrm{in}}\cup A)=\{\overline{x}^{(i)} : 0\leq i\leq k\}$. Since
$\Xi'_{\mathrm{in}}$ is $S_x(\Xi'_{\mathrm{in}})$-cartesian, by Theorem \ref{CDZ}, the
relation $H_0(\Xi'_{\mathrm{in}}\cup A)\supseteq H_1(\Xi'_{\mathrm{in}}\cup
A)\supseteq\cdots\supseteq H_{n_0}(\Xi'_{\mathrm{in}}\cup A)$ holds. From the description
of \textbf{MCS4}, we can deduce that $H_{n_0}(\Xi'_{\mathrm{in}}\cup A)\supseteq
H_{n_0+1}(\Xi'_{\mathrm{in}}\cup A)$, which leads to
\begin{eqnarray}\label{con1}
H_0(\Xi'_{\mathrm{in}}\cup A)\supseteq H_1(\Xi'_{\mathrm{in}}\cup
A)\supseteq\cdots\supseteq H_{n_0+1}(\Xi'_{\mathrm{in}}\cup A).
\end{eqnarray}
Note that for any $\overline{x}^{(i)}, 0\leq i\leq k$, there exists $h_i\in \{0, 1,
\ldots, m_{n_0}\}$ such that $\overline{x}^{(i)}=x_{h_i}$. Therefore, we could find a
permutation $\sigma$ of $\{0, 1, \ldots, m_0\}$ satisfying $\sigma(i)=h_i, i=0,\ldots,k,$
and $ \sigma(i)=i, i=m_{n_0}+1,\ldots,m_0$. Choose lines $l_i^y: x=x_{\sigma(i)},
i=0,\ldots,m_0$, that give rise to $V_i(\Xi'_{\mathrm{in}})=\{y_j : 0\leq j \leq
n_{\sigma(i)}\}, i=0,\ldots,m_0$. Since $n_0=n_1=\cdots=n_{m_{n_0}}$, the relation
$V_0(\Xi'_{\mathrm{in}})=V_1(\Xi'_{\mathrm{in}})=\cdots=V_{m_{n_0}}(\Xi'_{\mathrm{in}})\supseteq
V_{m_{n_0}+1}(\Xi'_{\mathrm{in}})\supseteq\cdots \supseteq V_{m_0}(\Xi'_{\mathrm{in}})$
holds. Observing that $V_i(\Xi'_{\mathrm{in}}\cup
A)=V_i(\Xi'_{\mathrm{in}})\cup\{\overline{y}\}, i=0, \ldots, k$, and
$V_i(\Xi'_{\mathrm{in}}\cup A)=V_i(\Xi'_{\mathrm{in}}), i=k+1,\ldots,m_0$, it is easy to
get
\begin{eqnarray*}
 V_0(\Xi'_{\mathrm{in}}\cup A)=\cdots=V_k(\Xi'_{\mathrm{in}}\cup
A)\supseteq V_{k+1}(\Xi'_{\mathrm{in}}\cup
A)\supseteq\cdots\supseteq V_{m_0}(\Xi'_{\mathrm{in}}\cup A).
\end{eqnarray*}
Thus together with (\ref{con1}), $\Xi'_{\mathrm{out}}=\Xi'_{\mathrm{in}}\cup A$ is
cartesian due to Theorem \ref{CDZ}, hence our statement is true.

For the $M$th loop, if $\Xi=[\ ]$, then $\Xi'$ here equals to the $\Xi'_{\mathrm{out}}$
of the \textbf{MCS4} of the $(M-1)$th loop that is cartesian due to the statement above.
Otherwise, since the algorithm stops in \textbf{MCS3} of this loop, $\Xi$ is a non-empty
cartesian set. Similar to the arguments above, we can prove that
$\Xi'=\Xi'_{\mathrm{out}}\cup \Xi$ is also cartesian.

Finally, we should verify that the output $\Xi'$ of the algorithm is maximal. Otherwise,
there must exist a maximal $S_x(\Xi'')$-cartesian subset $\Xi''$ of $\Xi$ satisfying
$\Xi''\supsetneqq\Xi'$. Take a point $\xi_0=(x_{i_0},y_{j_0})$ with
$(i_0,j_0)=\min_{\prec_{\mathrm{inlex}}}\{(i,j)\in S_x(\Xi'') :
(x_i,y_j)\in\Xi''\setminus\Xi'\}$. Suppose there exists a point in $\Xi'$ sharing the
ordinate with $\xi_0$. If it is chosen as a point in the maximal row subset in
\textbf{MCS4} of some loop, by the definition of $\xi_0$, we know that $\xi_0$ is surely
contained in the set $\Xi$ of that step, which contradicts the definition of the maximal
row subset. Otherwise, it must appear in the cartesian set $\Xi$ in \textbf{MCS3} in the
final loop. Then, by the definition of $\xi_0$, it should be contained in $\Xi$ hence the
output set $\Xi'$, which introduces a contradiction. If there does not exist a point in
$\Xi'$ sharing the ordinate with $\xi_0$, since $\Xi''$ is also cartesian, by Theorem
\ref{CDZ}, it is easily to see that $\xi_0$ must remain in $\Xi$ in every loop, which
contradicts the termination condition. As a result, the output of the Algorithm
\ref{MaxSubset} is a maximal cartesian subset. \qed

\end{pf}


Let us continue with the setup and notation in Algorithm \ref{MaxSubset}, and assume that
the final output of it is $\Xi'$ who is $S_x(\Xi')$-cartesian . We now discuss how to
preprocess the BM algorithm with the help of $\Xi'$.

Define an order $\prec_\Xi$ on the set $\Xi$. Let $\xi^{(1)}, \xi^{(2)}\in \Xi$. We say
that $\xi^{(1)}\prec_\Xi \xi^{(2)}$ if one of the following conditions holds:

\begin{itemize}
  \item [(1)]$\xi^{(1)}\in\Xi'$, and
  $\xi^{(2)}\in \Xi\backslash\Xi'$.
  \item [(2)]$\xi^{(1)}=(x_{i_1}, y_{j_1}), \xi^{(2)}=(x_{i_2},
y_{j_2})\in \Xi'$ and $(i_1, j_1)\prec_{\mathrm{inlex}}(i_2, j_2)$ with $(i_k, j_k)\in
S_x(\Xi'), k=1,2$.
\end{itemize}
It should be noticed that the order is not total. For the points in $\Xi\backslash\Xi'$,
any order of them can be interpreted as increasing. Hereafter, we will suppose that the
points in $\Xi=\{\xi^{(1)}, \ldots, \xi^{(\#\Xi)}\}$ have been ordered increasingly
w.r.t. $\prec_\Xi$, namely $\xi^{(i)}\prec_\Xi \xi^{(j)}, 0\leq i<j\leq \#\Xi$. By the
definition of $\prec_\Xi$, we have $\Xi'=\{\xi^{(1)}, \ldots, \xi^{(\#\Xi')}\}$.

According to Lemma \ref{subcon},  $N'=\{x^iy^j : (i, j)\in S_x(\Xi')\}\subset N$, with
$N$ as a member of the 3-tuple output of BM algorithm. Thus the other monomials of $N$
are obviously contained in $\mathbb{T}^2\backslash N'$. Notice that the generators of
$\mathbb{T}^2\backslash N'$ are located in the border of $N'$, denoted by $L$, we can
continue to spot the elements in $L$ by BM algorithm to complete $N$.

Next, we will pay attention to the computation of the Newton basis. Since $\Xi'$ is
cartesian,  recalling Proposition \ref{phi^x}, we can construct the polynomials
$\phi_{ij}^x$ w.r.t. $S_x(\Xi')$. Order $\phi_{ij}^x, (i, j)\in S_x(\Xi')$, increasingly
w.r.t. $(i, j)$ under $\prec_{\mathrm{inlex}}$, and denote them as $q_1, q_2, \ldots,
q_{\#\Xi'}$. Set the matrix
\begin{equation}\label{matrix}
  B=\left( \begin{array}{cccc}
    q_1(\xi^{(1)}) & q_1(\xi^{(2)}) & \cdots & q_1(\xi^{(\#\Xi')}) \\
    q_2(\xi^{(1)}) & q_2(\xi^{(2)}) & \cdots & q_2(\xi^{(\#\Xi')}) \\
    \vdots & \vdots &  & \vdots \\
    q_{\#\Xi'}(\xi^{(1)}) & q_{\#\Xi'}(\xi^{(2)}) & \cdots & q_{\#\Xi'}(\xi^{(\#\Xi')}) \\
  \end{array}
\right).
\end{equation}
By Proposition \ref{phi^x}, $B$ is obviously upper unitriangular which implies that the
polynomials $q_1, q_2, \ldots, q_{\#\Xi'}$ constitute a Newton basis for
$\mathcal{P}_{\prec}(\Xi')=\mathrm{Span}_{\mathbb{F}}N'$.

All in all, with the notation above, we get our preprocessing procedure for BM algorithm.

\begin{alg}\label{prBM}
(GPBM)

\textbf{Input}: A set of distinct points $\Xi\subset \mathbb{F}^2$ and a term order
$\prec$.

\textbf{Output}: The 3-tuple $(G,N,Q)$.

\vskip 3mm

\indent \textbf{GPBM1}: Get a maximal cartesian subset $\Xi'$ of $\Xi$ by Algorithm
\ref{MaxSubset};
\\
\indent \textbf{GPBM2}: Compute the lower set $S_x(\Xi')$ w.r.t.  $\Xi'$, the set
$N:=\{x^iy^j : (i,j)\in S_x(\Xi')\}$, and the set $Q:=\{q_1,q_2,\ldots,q_{\#\Xi'}\}$
where the $q_i$'s are as in (\ref{matrix}).
\\
\indent \textbf{GPBM3}: Construct $L:=\{x\cdot t : t\in N\}\bigcup\{ y\cdot t : t\in
N\}\setminus N$ and the matrix $B$ that is same to \eqref{matrix}.
\\
\indent \textbf{GPBM4}: Goto \textbf{BM2} of the BM algorithm to complete the computation
and get the whole output.
\end{alg}

\section{Implementation and Timings}\label{Application}

From the above section,  we can see easily that our preprocessing paradigm is more
suitable to the cases where the constructed maximal cartesian subset $\Xi'$ forms a
relatively large proposition in $\Xi$. Especially, when the field $\mathbb{F}$ is finite,
our preprocessing will play a more important role in consideration of the nature of
finite fields. In this section, we will present some experimental results to compare the
effectiveness of our paradigm with the classical BM.  First see an example with point set
of small size.

\begin{exmp}
We choose the field $\mathbb{F}_7$, and let
\begin{align*}
    \Xi=\{&(0, 0), (0, 1), (0, 4), (0, 5), (1, 0), \\
          &(1, 1), (1, 4), (1, 6), (2, 1), (2, 2), \\
          &(2, 6), (3, 2), (4, 2), (4, 5), (4, 6), \\
          &(5, 1), (5, 5), (5, 6), (6, 0), (6, 2)\}.
\end{align*}

By Algorithm \ref{MaxSubset}, we can construct the maximal cartesian subset
$$
\Xi'=\{(0, 1), (1, 1), (2, 1), (5, 1), (1, 6), (2, 6), (5, 6), (1, 0), (1, 4)\}
$$
hence get
\begin{align*}
N=&\{1, x, x^2, x^3, y, xy, x^2y, y^2\},\\
Q=&\{1, x, 4x^2+3x, 2x^3+x^2+4x, 3y+4, 3xy+4x+4y+3,\\
   &\phantom{\{}2x^2y+5x^2+xy+6x+4y+3, 6y^2+1, 2y^3+5y\},\\
L=&\{y^4, xy^2, xy^3, x^2y^2, x^3y, x^4\},\\
B=&\left(
     \begin{array}{cccc}
       1 & 1 & 1 & \cdots \\
       0 & 1 & 2 & \cdots \\
       0 & 0 & 1 & \cdots \\
       \vdots & \vdots & \vdots & \ddots \\
     \end{array}
   \right).
\end{align*}
Put these $N, Q, L, B$ into BM algorithm, we can get the final output
\begin{align*}
N=&\{1, x, x^2, x^3, y, xy, x^2y, y^2, y^3, xy^2, y^4, xy^3, x^2y^2, x^3y, x^4, y^5, xy^4, x^2y^3,\\
    &\phantom{\{} x^3y^2, x^4y\},\\
Q=&\{ 1, x, 4x^2+3x, 2x^3+x^2+4x, 3y+4, 3xy+4x+4y+3, \\
  &\phantom{\{}2x^2y+5x^2+xy+6x+4y+3, 6y^2+1, 2y^3+5y, xy^2+6y^2+6x+1,\\
  &\phantom{\{} y^4+3y^3+6y^2+4y, 5xy^3+5y^4+3y^3+2xy+2y^2+4y,\\
  &\phantom{\{}6x^2y^2+xy^2+x^2+6x,\ldots \},\\
\end{align*}
\begin{align*}
G=&\{y^6+3y^5+2y^4+6y^3+4y^2+5y,\\
  &\phantom{\{}xy^5+x^4y+6x^3y^2+x^2y^3+5xy^4+6y^5+6x^4+2x^3y+6x^2y^2+3xy^3+\\
  &\phantom{\{}3y^4+6x^3+6x^2y+2xy^2+6y^3+x^2+2xy+6y^2+x,\\
  &\phantom{\{}x^2y^4+x^4y+3x^2y^3+3xy^4+5y^5+x^4+6x^3y+3x^2y^2+2xy^3+4y^4+\\
  &\phantom{\{}6x^3+4y^3+6x^2+2xy+3y^2+x+5y,\ldots\}.\\
\end{align*}

\end{exmp}

In the following, several tables show the timings for the computations of BM-problems on
sets of distinct random points w.r.t. the term order $\prec_{\mathrm{lex}}$ or
$\prec_{\mathrm{tdinlex}}$. The algorithms presented in the paper were implemented on
Maple 12 installed on a laptop with 2 Gb RAM and 1.8 GHz CPU.

Take the field $\mathbb{F}_{23}$,  we have
\begin{center}
\begin{tabular}{lllll}\hline
  $\#\Xi$ & 200     & 300      & 400      & 500\\ \hline\hline
  BM      & 4.968 s & 15.359 s & 34.609 s & 61.172 s\\
  SPBM     & 1.438 s & 3.766 s  & 7.141 s  & 7.969 s\\
  \hline
\end{tabular}
\end{center}

For $\mathbb{F}_{37}$, we have
\begin{center}
\begin{tabular}{lllll}\hline
  $\#\Xi$          & 300      & 600       & 900       & 1200\\ \hline\hline
  BM               & 16.265 s & 121.766 s & 420.219 s & 1060.203 s\\
  SPBM              & 4.172 s  & 25.125 s  & 82.000 s  & 132.719 s\\
  \hline
\end{tabular}
\end{center}

For $\mathbb{F}_{17}$, we have
\begin{center}
\begin{tabular}{cllll}\hline
  $\#\Xi$    & 100     & 150     & 200     & 250\\ \hline\hline
  BM         & 0.875 s & 2.421 s & 4.953 s & 8.188 s\\
  GPBM        & 0.797 s & 2.125 s & 4.250 s & 5.641 s\\ \hline
  Preprocessing & 0.015 s & 0.094 s & 0.172 s & 0.391 s\\
  \hline
  $\#\Xi'/\#\Xi$      & 0.310   & 0.393   & 0.430   & 0.616 \\ \hline
\end{tabular}
\end{center}

Take the field $\mathbb{F}_{29}$, we have
\begin{center}
\begin{tabular}{cllll}\hline
  $\#\Xi$      & 200     & 400      & 600       & 800      \\ \hline\hline
  BM           & 5.672 s & 38.063 s & 112.156 s & 235.813 s \\
  GPBM          & 5.562 s & 36.906 s & 105.828 s & 135.609 s \\ \hline
  Preprocessing   & 0.046 s & 0.313 s  & 1.671 s   & 8.125 s\\
  \hline
  $\#\Xi'/\#\Xi$        & 0.125   & 0.178    & 0.328     & 0.711     \\ \hline
\end{tabular}
\end{center}

\bibliographystyle{elsarticle-num}
\bibliography{ref}

\begin{thebibliography}{10}
\expandafter\ifx\csname url\endcsname\relax
  \def\url#1{\texttt{#1}}\fi
\expandafter\ifx\csname urlprefix\endcsname\relax\def\urlprefix{URL }\fi
\expandafter\ifx\csname href\endcsname\relax
  \def\href#1#2{#2} \def\path#1{#1}\fi

\bibitem{CLO2005}
D.~A. Cox, J.~Little, D.~O'Shea, Using Algebraic Geometry, 2nd Edition, Vol.
  185 of Graduate Texts in Mathematics, Springer, New York, 2005.

\bibitem{Sau2006}
T.~Sauer, Polynomial interpolation in several variables: Lattices, differences,
  and ideals, in: K.~Jetter, M.~Buhmann, W.~Haussmann, R.~Schaback,
  J.~St\"{o}ckler (Eds.), Topics in Multivariate Approximation and
  Interpolation, Vol.~12 of Studies in Computational Mathematics, Elsevier,
  Amsterdam, 2006, pp. 191--230.

\bibitem{dBo2005}
C.~de~Boor, Ideal interpolation, in: C.~K. Chui, M.~Neamtu, L.~L. Schumaker
  (Eds.), Approximation Theory XI: Gatlinburg 2004, Nashboro Press, Brentwood
  TN, 2005, pp. 59--91.

\bibitem{Sak1998}
S.~Sakata, Gr\"{o}bner bases and coding theory, in: B.~Buchberger, F.~Winkler
  (Eds.), Gr\"{o}bner Bases and Applications, Vol. 251 of London Mathematical
  Society Lecture Notes Series, Cambridge University Press, New York, 1998, pp.
  470--485.

\bibitem{Sal2009}
M.~Sala, Gr\"{o}bner bases, coding, and cryptography: a guide to the
  state-of-art, in: M.~Sala, T.~Mora, L.~Perret, S.~Sakata, C.~Traverso (Eds.),
  Gr\"{o}bner Bases, Coding, and Cryptography, Springer, Berlin, 2009, pp.
  1--8.

\bibitem{Rob1998}
L.~Robbiano, Gr\"{o}bner bases and statistics, in: B.~Buchberger, F.~Winkler
  (Eds.), Gr\"{o}bner Bases and Applications, Vol. 251 of London Mathematical
  Society Lecture Notes Series, Cambridge University Press, New York, 1998, pp.
  179--204.

\bibitem{LS2004}
R.~Laubenbacher, B.~Stigler, A computational algebra approach to the reverse
  engineering of gene regulatory networks, J. Theoret. Biol. 229~(4) (2004)
  523--537.

\bibitem{JS2006}
W.~Just, B.~Stigler, Computing {Gr\"{o}bner} bases of ideals of few points in
  high dimensions, ACM Commun. Comput. Algebra 40~(3-4) (2006) 67--78.

\bibitem{MB1982}
H.~M\"{o}ller, B.~Buchberger, The construction of multivariate polynomials with
  preassigned zeros, in: J.~Calmet (Ed.), Computer Algebra: EUROCAM '82, Vol.
  144 of Lecture Notes in Computer Science, Springer, Berlin, 1982, pp. 24--31.

\bibitem{FGLM1993}
J.~C. Faug\`{e}re, P.~Gianni, D.~Lazard, T.~Mora, Efficient computation of
  zero-dimensional {Gr\"{o}bner} bases by change of ordering, J. Symbolic
  Comput. 16~(4) (1993) 329--344.

\bibitem{MMM1993}
M.~G. Marinari, H.~M. M\"{o}ller, T.~Mora, Gr\"{o}bner bases of ideals defined
  by functionals with an application to ideals of projective points, Appl.
  Algebra Engrg. Comm. Comput. 4~(2) (1993) 103--145.

\bibitem{ABKR2000}
J.Abbott, A.Bigatti, M.Kreuzer, L.Robbiano, Computing ideals of points, J.
  Symbolic Comput. 30 (2000) 341--356.

\bibitem{CM1995}
L.~Cerlienco, M.~Mureddu, From algebraic sets to monomial linear bases by means
  of combinatorial algorithms, Discrete Math. 139~(1-3) (1995) 73--87.

\bibitem{GRS2003}
S.~Gao, V.~M. Rodrigues, J.~Stroomer, Gr\"{o}bner basis structure of finite
  sets of points (Preprint).

\bibitem{FRR2006}
B.~Felszeghy, B.~R\'{a}th, L.~R\'{o}nyai, The lex game and some applications,
  J. Symbolic Comput. 41~(6) (2006) 663--681.

\bibitem{Sau2004}
T.~Sauer, Lagrange interpolation on subgrids of tensor product grids, Math.
  Comp. 73~(245) (2004) 181--190.

\bibitem{Cra2004}
N.~Crainic, Multivariate {Birkhoff-Lagrange} interpolation schemes and
  cartesian sets of nodes, Acta Math. Univ. Comenian.(N.S.) LXXIII~(2) (2004)
  217--221.

\bibitem{CDZ2006}
T.~Chen, T.~Dong, S.~Zhang, The {Newton} interpolation bases on lower sets, J.
  Inf. Comput. Sci. 3~(3) (2006) 385--394.

\bibitem{BW1993}
T.~Becker, V.~Weispfenning, Gr\"{o}bner Bases, Vol. 141 of Graduate Texts in
  Mathematics, Springer-Verlag, New York, 1993.

\bibitem{Lor1992}
R.~Lorentz, Multivariate Birkhoff Interpolation, Vol. 1516 of Lecture Notes in
  Mathematics, Springer, Heidelberg, 1992.

\bibitem{dBo2007}
C.~de~Boor, Interpolation from spaces spanned by monomials, Adv. Comput. Math.
  26~(1) (2007) 63--70.

\end{thebibliography}

\end{document}